\documentclass[sts]{imsart}

\RequirePackage[OT1]{fontenc}
\usepackage{amsthm,amsmath,amsfonts,amssymb,natbib}
\RequirePackage[colorlinks,citecolor=blue,urlcolor=blue]{hyperref}

\startlocaldefs
\numberwithin{equation}{section}
\theoremstyle{plain}
\newtheorem{thm}{Theorem}[section]
\newtheorem{lemma}{Lemma}[section]
\newtheorem{example}{Example}[section]

\endlocaldefs

\newcommand{\btheta}{\mbox{\boldmath $\theta$}}

\newcommand{\R}{\mbox{$\mathbb{R}$}}
\newcommand{\GG}{\mbox{$\mathbb{G}$}}

\newcommand{\dis}{\mbox{$D_{\sc KW}$}}

\newcommand{\pr}{\mbox{\sc pr}}

\newcommand{\sumin}{\sum_{i=1}^n}
\newcommand{\bea}{\begin{eqnarray}}
\newcommand{\eea}{\end{eqnarray}}
\newcommand{\ba}{\begin{eqnarray*}}
\newcommand{\ea}{\end{eqnarray*}}
\newcommand{\be}{\begin{equation}}
\newcommand{\ee}{\end{equation}}
\newcommand{\bi}{\begin{itemize}}
\newcommand{\ei}{\end{itemize}}

\newcommand{\var}{\mbox{\sc var}}

\newcommand{\EE}{\mbox{\rm E}}
\newcommand{\iid}{\mbox{i.i.d.\,}}
\newcommand{\cdf}{\mbox{c.d.f.\! }}

\newcommand{\no}{\noindent}
\newcommand{\vs}{\vspace{1ex}}

\begin{document}

\begin{frontmatter}
\title{Consistency of the MLE under mixture models
%\thanksref{T1}
}
\runtitle{Consistency of MLE}
%\thankstext{T1}{To be added if needed}

\begin{aug}
\author{\fnms{Jiahua} \snm{Chen}
%\thanksref{t1}
\ead[label=e1]{jhchen@stat.ubc.ca}}
%\author{\fnms{Second} \snm{Author}\ead[label=e2]{second@somewhere.com}}
%\and
%\author{\fnms{Third} \snm{Author}\thanksref{t1}
%\ead[label=e3]{third@somewhere.com}
%\ead[label=u1,url]{www.foo.com}}

%\thankstext{t1}{Some comment}
%\thankstext{t2}{First supporter of the project}
\runauthor{Chen}

\affiliation{Department of Statistics, Yunnan University, and University of British Columbia}

%\address{Address of the First and Second authors,
%usually few lines long \printead{e1,e2}.}

\end{aug}

\begin{abstract}
The large-sample properties of likelihood-based statistical inference
under mixture models have received much attention from statisticians.
Although the consistency of the nonparametric MLE is regarded as
a standard conclusion, many researchers ignore the precise conditions required on the 
mixture model. An incorrect claim of consistency can lead to
false conclusions even if the mixture model under investigation
seems well behaved.
Under a finite normal mixture model, for instance,
the consistency of the plain MLE is often erroneously assumed
in spite of recent research breakthroughs. 
This paper streamlines the consistency
results for the nonparametric MLE in general, and in particular for
the penalized MLE under finite normal mixture models.
\end{abstract}


\end{frontmatter}

\section{Introduction}

A family of distributions, parametric or nonparametric, is 
regarded as a probability or statistical model.
A parametric model can hence be a family of density functions 
in the form $\{f(x; \theta): \theta \in \Theta\}$
where $\Theta$ is a subset of ${\R}^d$ for some positive integer $d$.
The measure with respect to which these densities are defined
is $\mu$, but this will be de-emphasized.
In applications, statisticians must 
select one or a set of these distributions to characterize the uncertainty
displayed in a random sample from a population.

A specific distribution family or model is often tentatively 
chosen in accordance with the scientific background of the application. 
For instance, the Poisson distribution/model is a textbook choice for the
number of annual car accidents of a policy holder, but it may not be ideal 
because risk levels vary.
It can therefore be helpful to divide the population
into several subpopulations, each modeled 
by its own Poisson distribution.
Thus, a finite mixture of Poisson distributions is 
a better choice for the pooled accident data. 
This leads to a generic finite mixture model in the form
\begin{equation}
\label{mixture}
f(x; G) = \sum_{j=1}^m \alpha_j f(x; \theta_j).
\end{equation}
In this formulation,
$f(x; G)$ is a finite convex combination of the component density functions
$f(x; \theta)$. 
Each $f(x; \theta_j)$ is the density for a subpopulation.
The recipe for the mixture is summarized by $G$,
specified by its cumulative distribution function (\cdf)
\be
\label{mixing}
G(\theta) = \sum_{j=1}^m \alpha_j I(\theta_j \leq \theta)
\ee
where $I(\cdot)$ is an indicator function.
The model $f(x; G)$ in the above definition has order $m$, 
even if some $\theta_j$ are equal or some $\alpha_j = 0$. 
A finite mixture model of order $m$ is hence also a 
degenerate order $m+1$ mixture.

For any \cdf  $G$ on $\Theta$, the following integral
\begin{equation}
\label{mixtureNP}
f(x; G) = \int f(x; \theta) dG(\theta)
\end{equation}
is a well-defined density function. 
When $G$ has only a finite number of
support points, as in \eqref{mixing}, $f(x; G)$ reduces to the finite
mixture model \eqref{mixture}.
When the form of $G$ is unspecified, 
\eqref{mixtureNP} is regarded as a nonparametric
mixture model.
The space of all the mixing distributions on $\Theta$ is
denoted $\GG$, and the space of those 
with at most $m$ support points is $\GG_m$.
It can be seen that $\GG_m \subset \GG_{m+1}$.
See \cite{r5},  \cite{r8}, and \cite{r7}  for the
general theory and applications of mixture models.

Research into mixture models has a long history. 
The most cited early publication is \cite{r11};
%Pearson (1894) 
he used a two-component normal mixture model
for a biometric data set.
This model was motivated by
the apparent skewness in the ratio of forehead-width to body-length
of 1000 crabs from Naples.
It was suspected that the population contained two unidentified species. 
\cite{r11} employed the method of moments for the parameter estimation
and provided a detailed account of the numerical calculation 
in the absence of the modern computer.
Nowadays, we would estimate the parameters by the maximum likelihood estimator (MLE).
The MLE is favored
for its asymptotic efficiency under regular parametric models.
The unpleasant numerical obstacle is now a history.

Mixture distributions form a distinct class of nonregular statistical
models. 
They are notorious for presenting statisticians with serious technical challenges
as well as many pleasant surprises.
Given an \iid sample of size $n$ from a mixture model,
\cite{rr7} showed that the likelihood
ratio test statistic for homogeneity is stochastically unbounded 
defying the classical chisquared limiting distribution.
\cite{rr9} showed that the optimal rate for estimating
the mixing distribution is $O_p(n^{-1/4})$ as compared
with $O_p(n^{-1/2})$ for parameter estimation under
regular models.
\cite{rr6} interpreted some of these abnormal results 
by the degenerated Fisher information when the order 
of the finite mixture model is unknown.
At the same time,
\cite{r14} showed that the nonparametric MLE of $G$
has at most $k$ support points, the number of distinct observed values.
Given the order of the finite mixture model and an initial
mixing distribution $G^{(0)}$,  \cite{r15} found
that the famous EM-algorithm
leads to a sequence of $G^{(k)}$ that
converges to a local maximum of the likelihood
function and it can be a locally consistent MLE \cite{rr3}. 
Both \cite{r14} on the geometric property
and \cite{r15} on the algorithmic convergence
place only nominal conditions on $f(x; G)$.

There is also good news about the large-sample
properties of the MLE.
\cite{r1} and \cite{r2} proved that the nonparametric MLE of 
$G$ is strongly consistent in an \iid setting. 
Using quotient topology, \cite{r3}
gave another consistency proof under finite mixture models.

These consistency results are encouraging, but
careful discussions of the relationships
between these proofs are lacking. 
Why are there multiple proofs for a single consistency result?
How do these proofs relate to each other? 
Without a full understanding, researchers may cite a paper 
when its specific proof/result is not applicable to the target problem.
This paper streamlines the consistency proofs of the 
nonparametric MLE under mixture models and
crystallizes their applicability.
Another topic is consistent estimation under the 
finite normal mixture model: the consistency of the MLE under this
model is often erroneously assumed.
This paper aims to popularize the consistent
estimator based on the penalized likelihood.
The ultimate goal of this exercise is to 
provide a solid basis for research into the
large-sample properties of mixture models.

The remainder of the paper consists of two major sections.
Section~2 addresses classical consistency results for the
nonparametric MLE under a mixture model.
Section~3 addresses results related to the consistency
of the penalized (or regularized) MLE under the
finite normal mixture model. Section~4 provides concluding remarks.

\section{Classical consistency results}
Let  $x_1, \ldots, x_n$ be a set of \iid observations 
of size $n$ from a nonparametric mixture model $\{ f(x; G): G \in \GG \}$.
The likelihood function of $G$ is given by
\[
L_n(G) = \prod_{i=1}^n f(x_i; G).
\]
The log-likelihood function of $G$ is hence
\be
\label{eqn20.like}
\ell_n(G) = \sum_{i=1}^n \log f(x_i; G).
\ee
Both $L_n(G)$ and $\ell_n(G)$ are functions defined on $\GG$. 
The nonparametric MLE $\hat G_n$ of $G$ is a \cdf on 
$\Theta$ such that
\be
\label{eqn20.mle}
\ell_n(\hat G) = \sup \{ \ell_n(G): G \in \GG\}.
\ee
Most rigorously, $\hat G$ is related to
the sample in a measurable fashion,  and it
is potentially one of many possible global maximum points of
$\ell_n(\cdot)$. In addition, $\hat G$ can be a limiting
point of a sequence of mixing distributions $G_j$ such that 
$\ell_n(G_j) \to \sup \{ \ell_n(G): G \in \GG\}$
as $j \to \infty$. 

For ease of presentation, \eqref{eqn20.mle} 
is regarded as an unambiguous definition of MLE.
When $G$ is confined to $\GG_m$, this 
becomes the MLE of $G$ under finite mixture models.
When $m=1$ in $\GG_m$,
\eqref{eqn20.mle} defines the ordinary parametric MLE of
$\theta$ under the model $\{f(x; \theta): \theta \in \Theta\}$.

As acknowledged by \cite{r1}, their proof of
the consistency of $\hat G$ is a technical clone of that of 
\cite{r4}. One key technical preparation for this
proof is the following well-known inequality.

\begin{lemma} \mbox{\rm (Jensen's inequality)}.
\label{lemma20.1}
Let $X$ be a random variable such that $\mbox{\EE} |X| < \infty$
and let $\varphi(t)$ be a convex function. Then
\[
\EE\{ \varphi(X) \} \geq \varphi(\EE(X)).
\]
\end{lemma}

Suppose $X$ is a random variable with density
function $f(x; \theta)$ in general and $\theta^*$ is
the true value of the parameter.
Let $Y =f(X; \theta)/f(X; \theta^*)$.
It can be seen that
\[
\EE^*(Y) = \int \{ f(x; \theta) /f(x; \theta^*)\} f(x; \theta^*) d\mu(x) 
\leq \int f(x; \theta) d\mu(x) = 1,
\]
where $\EE^*$ is the expectation with respect to the distribution
with the specific $\theta^*$ parameter value.
Applying Jensen's inequality to $Y$ and 
$\varphi(t) = - \log (t)$, we get
\be
\label{JensenIn}
- \EE^* 
\log 
\{ f(X; \theta) /f(X; \theta^*)\} 
\geq - \log [ \EE^* \{ f(X; \theta) /f(X; \theta^*)\} ] 
\geq 0.
\ee
If $\EE^* | \log f(X; \theta^*)| < \infty$, then
\[
\EE^*\{ \log f(X; \theta)\} \leq \EE^* \{ \log f(X; \theta^*)\}.
\]
The equality holds if and only if $f(x; \theta) \equiv f(x; \theta^*)$,
or it is equal except for a zero-probability set of $x$ in terms of $f(x; \theta^*)$.
The expectation $\EE^*[ \log \{f(X; \theta^*)/f(X; \theta)\}]$
%%% changed G to \theta
is referred to as the Kullback--Leibler information 
between the two distributions.

The following lemma may seem trivial, but it is the
basis of most proofs for generic consistency.
%%% all changed to most.

\begin{lemma}
\label{lemma.trivial}
\mbox{\rm (Trivial Consistency of MLE)}.
Suppose the model under investigation is
$\{f(x; \theta), \theta \in \Theta\}$
where $\Theta = \{\theta^*, \theta_1, \ldots, \theta_M\}$.
In addition, $f(x; \theta) = f(x; \theta^*)$ except for
a zero-probability set of $x$ with respect to $f(x; \theta^*)$
implies $\theta = \theta^*$. 

Then the MLE $\hat \theta_n$ of $\theta$ equals $\theta^*$ almost surely
as $n \to \infty$.
\end{lemma}

\vs
\no
{\sc Proof}:
Let $X$ be a random variable from $f(x; \theta^*)$.
Then
\[
\max_{1 \leq j \leq M} \EE^* \log \{ f(X; \theta_j)/f(X; \theta^*)\} < 0.
\]
By the strong law of large numbers, this inequality implies
\be
\label{simpleMLE}
\max_{1 \leq j \leq M} \sumin \log f(x_i; \theta_j)
< 
 \sumin \log f(x_i; \theta^*)
 \ee
 almost surely.
By the definition of the MLE, 
$\hat \theta_n=\theta^*$ almost surely as $n \to \infty$.
\qed \vs

Surprisingly, most generic MLE consistency proofs are variations
or novel upgraded versions of this lemma.
%% all ---> most of??
These include the situations where $\Theta$ is a
subset of the Euclidean space $\R^d$, 
the space of all mixing distributions $\GG$ or $\GG_m$,
or any abstract space.

In Section \ref{sub2.1}, we first present a trivialized
Wald theorem, Theorem \ref{Wald}.
It is universal consistency result though not directly applicable to 
many practical models. However,
the full Wald Theorem can be easily explained from this angle.

Theorem \ref{Wald} is then used as basis to understand
and prove the consistency results of 
\cite{r1,r2} in the context of mixture models
subsequently presented in Sections \ref{sub2.2}--\ref{sub2.4}.

\subsection{Essentials of the Wald consistency proof}
\label{sub2.1}
Let $\{f(x; \theta), \theta \in \Theta\}$ be the probability model
under investigation, where $\Theta$ is an abstract but metric parameter
space. For any subset $B$ of $\Theta$, define
\[
f(x; B) = \sup_{\theta \in B} f(x; \theta).
\]
Denote the distance on $\Theta$ as $\rho(\cdot, \cdot)$.
For any $\epsilon > 0$, let
\[
B_\epsilon(\theta^*) = \{\theta:  \rho(\theta, \theta^*) < \epsilon \}
\]
be an open ball of radius $\epsilon$ centred at $\theta^*$. 
Let $B^c$ be the complement of $B$.

\begin{thm} 
\mbox{\rm (Trivialized Wald Theorem)}.
\label{Wald}
Let $x_1, \ldots, x_n$ be an \iid sample from 
$f(x; \theta^*)$, a member of $\{f(x; \theta),  \theta \in \Theta\}$.
Let $\hat \theta_n$ be an MLE of $\theta$
as defined in \eqref{eqn20.mle}.

Suppose that for any $\epsilon > 0$, there exists a finite number
of subsets  $B_1, B_2, \ldots, B_J$ of $\Theta$  such that
$B^c_\epsilon(\theta^*) \subset \{ \cup_{j=1}^J B_j\}$ and for each $j$
\be
\label{JensenG}
\EE^* \log \{ f(X; B_j)/f(X; \theta^*)\} < 0.  %% added /
\ee
Then, $\hat \theta_n \to \theta^*$ almost surely
as $n \to \infty$.
\end{thm}

The proof is identical to that of Lemma \ref{lemma.trivial}.
The technicalities of \cite{r4} amount to specifying conditions on the model 
that lead to the {\it generalized Jensen's inequality}
\eqref{JensenG}.
What are the conditions placed on $f(x; \theta)$ by \cite{r4}? 
Here is a reorganized and slightly altered list.

\bi
\item[(W1)]
\mbox{\rm Identifiability}: 
{\it Let $F(x; \theta)$ be the cumulative distribution function
of $f(x; \theta)$. $F(x; \theta) = F(x; \theta^*)$ for all $x$ implies $\theta = \theta^*$.
}
\ei

Clearly, if $F(x; \theta) = F(x; \theta^*)$ for all $x$, there 
will be no stochastic difference between the data sets generated 
from these two distributions. 
Hence, consistent estimation of $\theta$ based only on
data is not possible if $\theta \neq \theta^*$.
Identifiability is therefore a necessary condition
for the consistent estimation of $\theta$.  

\bi
\item[(W2)]
\mbox{\rm Continuity and slightly more}: 
{\it For all $x$,
$\lim_{\theta \to \theta_0} f(x; \theta) = f(x; \theta_0)$ for any given $\theta_0$
and $\lim_{|\theta| \to \infty} f(x; \theta) = 0$.
}
\ei

Here, $| \theta |$ can be the Euclidean norm or any norm appropriate in
the context of the problem.
Technically, continuity may not be required for a zero-probability 
set of $x$ in terms of $f(x; \theta_0)$. 
The models used in applications are usually well behaved.
Hence, this extra generality is not generally needed if
one takes appropriate care and therefore is not
included as part of the condition.

For the consistent estimation of a parameter,
the distributions with close $\theta$ values must be similar.
Therefore, the continuity condition is indispensable. 

\bi
\item[(W3)]
\mbox{\rm Finite Kullback--Leibler Information}:
Let $[\cdot]^+$ denote the positive part of the quantity in the brackets.
{\it
For any $\theta \neq \theta^*$, there exists an $\epsilon > 0$ such that
\[
\EE^* [ \log \{ f(X; B_\epsilon(\theta)) / f(X; \theta^*) \} ]^+< \infty
\]
 and there exists a large enough $r > 0$ such that
\[
\EE^*[ \log \{ f(X; B^c_r(\theta^*)) / f(X; \theta^*) \}]^+ < \infty.
\]
}
\ei

\bi
\item[(W4)]
\mbox{\rm Closeness}: 
{\it The parameter space $\Theta$ is a closed subset
of $\R^d$.}
\ei

The use of these conditions is demonstrated in the following lemma.

\begin{lemma}
\label{lemma2.3}
The conditions of Theorem \ref{Wald} are satisfied under (W1)--(W4).
\end{lemma}

\vs
\noindent
{\sc Proof}.
Note that $ \log \{ f(X; B_\epsilon(\theta)) / f(X; \theta^*) \}$ is a monotone
increasing function of $\epsilon$. 
The continuity condition (W2) ensures that 
$\lim f(X; B_\epsilon(\theta)) = f(X; \theta)$
as $\epsilon \to 0^+$, i.e., when $\epsilon$ decreases to 0. 
Thus, this condition validates the dominated convergence theorem
in the following way:
\[
\lim_{\epsilon \to 0^+}
\EE^* [ \log \{ f(X; B_\epsilon(\theta)) / f(X; \theta^*) \} ]^+
=
\EE^* [ \log \{ f(X;  \theta) / f(X; \theta^*) \} ]^+.
\]
For the negative counterpart of this expectation,
the famous Fatou's Lemma implies, under the continuity condition (W2),
\[
{\lim\inf}_{\epsilon \to 0^+} \EE^* [ \log \{ f(X; B_\epsilon(\theta)) / f(X; \theta^*) \} ]^- 
\geq
\EE^* [ \log \{ f(X; \theta) / f(X; \theta^*) \} ]^-.
\]
The monotonicity on the left-hand side in $\epsilon$ ensures that the
limit exists, rather than merely $\lim\inf$.
Hence,
\[
\lim_{\epsilon \to 0^+}
\EE^* [ \log \{ f(X; B_\epsilon(\theta)) / f(X; \theta^*) \} ]
\leq
\EE^* [ \log \{ f(X;  \theta) / f(X; \theta^*) \} ] <0
\]
with the strict $< 0$ implied by the identifiability condition (W1).

The above result shows that for each $\theta \neq \theta^*$, 
there is a small enough $\epsilon_\theta$ such that
\[
\EE^* \log \{ f(X; B_{\epsilon_\theta} (\theta)) / f(X; \theta^*) \}  < 0.
\]

Since $f(x; \theta) \to 0$ as $|\theta| \to \infty$, we can similarly show
that there exists $r$ such that
\[
\EE^* \log \{ f(X; B^c_r(\theta^*)) / f(X; \theta^*) \} < 0.
\]

For any $\epsilon > 0$, let $\bar{B}_r(\theta^*)$ be the closure of
${B}_r(\theta^*)$. Then
\[
A = \{ B^c_\epsilon(\theta^*)\}  \cap \{\bar {B}_r(\theta^*) \}
\]
is bounded and closed and therefore compact.
It is trivial to see that
\[
\big [ \cup_{\theta \in A} \{ B_{\epsilon_\theta} (\theta) \} \big ] \supset A.
\]
When a compact set is covered by the union of a collection of open sets,
it is covered by a finite number of such sets.
Applying this result here, we have a finite
number of $B_{\epsilon_\theta} (\theta)$ whose union covers
the compact set $A$. Let these be $B_1, \ldots, B_{J-1}$,
and reserve $B_J$ for $B_r(\theta^*)$.

All the conditions in Theorem \ref{Wald} are satisfied.
Hence, the MLE is consistent under (W1)--(W4). 
\qed\vs

We have successfully established the consistency result
of \citep{r4} stated as Lemma \ref{lemma2.3}.
The proof is markedly simplified and it is obtained without requiring
\[
E^* | \log f(X; \theta^*)| < \infty.
\]

\vs\vs
When $\Theta$ is not a closed set, it is often possible to 
continuously extend
the range of the density function $f(x; \theta)$
to the closure of $\Theta$,  $\bar{\Theta}$. 
The Wald Theorem may then be applicable to the expanded
model.
Take the Poisson distribution family as an example:
\[
f(x; \theta) = \frac{\theta^x}{x!} \exp( - \theta) 
\]
for $x=0, 1, \ldots$ and $\Theta = (0, \infty)$.
It can be seen that $\Theta$ is not a closed set. However, by defining
$f(x; 0) = 0$ for all $x$ except for $f(x; 0) = 1$ when $x =0$,
we extend the model to $\bar{\Theta} = [0, \infty)$. 
Consequently, the conditions for a consistent MLE are satisfied
and the MLE is consistent.
Of course, the consistency of the MLE can be
easily established without utilizing the generic
Wald Theorem.

\subsection{Consistency of nonparametric MLE under mixture model: 
Kiefer--Wolfowitz (KW) approach}
\label{sub2.2}

%Under a nonparametric mixture model, the density
%function under investigation is $f(x; G)$ with infinite-dimensional 
%parameter space $\GG$. I wish to demonstrate that
%Theorem \ref{Wald} is applicable after appropriately
%interpreting (W1)--(W4) in the new context.
%This is all \cite{r1} is about in a nutshell.

We now illustrate and understand the contribution of \cite{r1}
through its connection to Theorem \ref{Wald}.
To discuss the consistence of the nonparametric MLE,
we need to choose a distance $D(\cdot, \cdot)$ on $\GG$
such as
\[
\dis(G_1, G_2) = \int_\Theta |G_1(\theta) - G_2(\theta)| \exp (- |\theta|) d \theta
\]
where $|\theta|$ is interpreted as $|\theta_1| + \cdots + |\theta_d|$
and $d \theta$ as $d\theta_1 \cdots d\theta_d$ when $\Theta \subset \R^d$.
Note that there are two $d$'s involved:
one is the dimension of $\theta$ and 
the other is the differential symbol for the integration. 

At this moment, a generic notion of distance suffices.
We say that $G \to G_0$ if $D(G, G_0) \to 0$.
Suppose $G^* \in \GG$ is the true mixing distribution
and $\hat G$ is an estimator. Then
$\hat G$ is strongly consistent when $D(\hat G, G^*) \to 0$
almost surely.

\bi
\item[(KW1)]
\mbox{\rm Identifiability}: 
{\it Let $F(x; G)$ be the cumulative distribution function
of $f(x; G)$. If $F(x; G) = F(x; G^*)$ for all $x$, then $D(G, G^*)=0$.
}
\ei

Suppose $\dis(\cdot, \cdot)$ is chosen as the distance on $\GG$.
When $f(x; \theta)$ is confined to the Poisson distribution, (KW1)
is satisfied. When $f(x; \theta)$ is the normal distribution, this condition
is violated. When $f(x; \theta)$ is binomial, it is also violated
in general. When $\GG$ is reduced to $\GG_m$, the normal mixture
satisfies (KW1); the binomial mixture
satisfies (KW1) when $m$ is small. See \cite{rr1,rr2}.

For $\GG_m$,
a mixing distribution can be expressed by two vectors: 
one for the component parameter values 
and one for the corresponding mixing proportions.
Suppose the Euclidean distance on this vector
space is chosen.
Then $f(x; G)$ remains the same 
when the entries of the two vectors are permuted.
This is the loss of identifiability due to label switching.
Quotient topology as suggested by \cite{r3} can be used to
avoid this problem. 
Label switching leads to technical difficulties for Bayesian analysis;
see \cite{r20} and \cite{r21}.

\bi
\item[(KW2)]
\mbox{\rm Continuity}: 
{\it The component parameter space $\Theta$ is a closed set.
For all $x$ and any given $G_0$,
\[
\lim_{G \to G_0} f(x; G) = f(x; G_0).
\]
}
%%% G is replace by \theta.

\item[(KW3)]
\mbox{\rm Finite Kullback--Leibler Information}:
{\it
For any $G \neq G^*$, there exists an $\epsilon > 0$ such that
\[
\EE^*[ \log \{ f(X; B_\epsilon(G)) / f(X; G^*) \}]^+ < \infty.
\]
}

\item[(KW4)]
\mbox{\rm Compactness}: {\it The definition of the mixture density
$f(x; G)$ in $G$ can be continuously extended to a compact
space $\bar {\GG}$ while retaining the validity of (KW3).}
\ei
%
%Unlike the corresponding condition under a parametric model,
%$B_\epsilon(G)$ in (KW3) needs some explanation.
%Superficially, it is a set of mixing distributions within 
%$\epsilon$-distance of $G$. This leads to the need for a distance
%on $\GG$ and the compact $\bar{\GG}$. 
%Let us assume that such a distance is available; we state the
%key result first.

\begin{thm}
\label{KWthm}
Let $x_1, \ldots, x_n$ be an \iid sample from 
$f(x; G^*)$, a member of $\{f(x; G), G\in \GG\}$.
Suppose conditions (KW1)--(KW4) are satisfied.
Then, the nonparametric MLE $\hat G_n$ is strongly consistent.
Namely, $\dis(\hat G, G^*) \to 0$ almost surely.
\end{thm}

\noindent
{\sc Proof}. Under the theorem conditions and following the
proof of Theorem \ref{Wald}, there must exist
a $\delta > 0$ for each $G \neq G^*$ such that
\[
\EE^*[ \log \{ f(X; B_\delta(G)) / f(X; G^*) \}] < 0.
\]
This leads to a finite open cover of the compact
set $B^c_\epsilon(G^*)$ for any $\epsilon > 0$.
Hence, by the law of large numbers,
\[
\max _{G \not \in B_\epsilon(G^*)} \ell_n(G) 
\leq
\ell_n(G^*) 
\]
almost surely. The arbitrariness of $\epsilon$ implies that
$\hat G_n$ is within an infinitesimal neighborhood of $G^*$
almost surely as $n \to \infty$ and is therefore consistent for
$G^*$.
\qed \vs

The proof of the consistency result in Theorem \ref{KWthm} 
is not fundamentally different from that of \cite{r1}.
The simplicity comes from requiring high level conditions
(KW3) and (KW4), and from not repeating some steps
in Lemma \ref{lemma2.3}. 
The current proof promotes the understanding
of the essence of their proof.
At the same time, (KW3) and (KW4) can be established
by (a): introducing a specific distance on $\GG$; (b): extending
the domain of $f(x; G)$ in $G$ continuously to a compact $\bar {\GG}$;
(c): verifying that Jensen's inequality hold on $\bar {\GG}$.
Going over these steps brings back the complexity.
We illustrate this point subsequently based on $\dis$.
It is known that $G_m \to G$ in distribution if and only
if $\dis(G_m, G) \to 0$ as $m \to \infty$.

\subsubsection{Compactificaton of $\GG$.}

Based on $\dis(\cdot, \cdot)$, the distance is no more than 
$\int_\Theta  \exp (- |\theta|) d \theta < \infty$.
To compactify $\GG$, we introduce
$\rho G$ as a subdistribution for any $\rho \in [0, 1)$. 
Let $\bar{\GG}$ be $\GG$ supplemented with all the subdistributions. 
Clearly, $\dis(\cdot, \cdot)$ is easily extended to $\bar{\GG}$.
Similarly, we naturally extend the range of $f(x; G)$ 
to $G \in \bar{\GG}$ by defining
\[
f(x; \rho G) = \rho \int_\Theta f(x; \theta) dG(\theta).
\]
Although we have technically expanded $\GG$ to $\bar{\GG}$, the
likelihood cannot be maximized on $\bar{\GG} - \GG$ since if
$f(x; G) \neq 0$ then $f(x; \rho G) < f(x; G)$ when $\rho < 1$. 
Hence, the MLE on
$\bar{\GG}$ is always a proper distribution on $\Theta$.

Since $\Theta$ is a closed subset of $\R^d$ from condition (KW2),
the limit of any Cauchy sequence in $\GG$
in terms of $\dis(\cdot, \cdot)$ is a subdistribution on $\Theta$.
Hence, $\bar{\GG}$ is a closed set. With the addition of total boundedness,
$\bar{\GG}$ is compact.

The extension from $\{f(x; G): G \in \GG\}$ to $\{f(x; G): G \in \bar{\GG}\}$
is largely symbolic. The real issues are continuity and (KW3)
on $\bar{\GG}$. For instance, the extension would fail on a normal
mixture because the component density function is not defined
at $\sigma = 0$. Otherwise, the continuous
extension is feasible.

\subsubsection{Continuity of $f(x; G)$ on $\bar{\GG}.$} 
Here is the sufficient and likely also necessary condition for
the continuity of $f(x; G)$ on $\bar{\GG}$ based on $\dis$
distance. Recall that $\dis(G_m, G) \to 0$ if and
only if $G_m \to G$ in distribution.

\begin{lemma}
Under (W2) and (W4), the extended mixture model
$\{ f(x; G): G \in \bar{\GG}\}$ is continuous in $G$ for all given $x$.
\end{lemma}

\no
{\sc Proof}: 
Recall that $G_m \to G_0$ in distribution if and only if
$\int h(\theta) d G_m(\theta) \to \int h(\theta) dG_0(\theta)$
for all bounded and continuous functions $h(\cdot)$,
according to one of many equivalent definitions.
By condition (W2),
$\lim_{|\theta| \to \infty} f(x; \theta) = 0$. Thus, for given $x$,
$f(x; \theta)$ is continuous and bounded on $\Theta$.
Therefore, this definition leads to
\[
f(x; G_m) = \int f(x; \theta) dG_m(\theta)
\to 
 \int f(x; \theta) dG_0(\theta) = f(x; G_0) 
\]
whenever $G_m \to G_0$ in distribution for all $G_0 \in \bar{\GG}$.
\qed\vs

Remark: Here the convergence in distribution 
includes subdistributions.
%%%

\subsubsection{Generalized Jensen's inequality.}
Technicality is inevitable when it comes to 
establishing the generalized Jensen's inequality 
for the mixture model on $\bar{\GG}$.
The user must decide whether or not
(KW3) holds on $\bar{\GG}$ for each specific mixture model. 
Two examples will help to explain the issue.

\begin{example}
If $G^*(M) = 1$ for some $M < \infty$, then the Poisson
mixture model satisfies all the conditions of Theorem \ref{KWthm}.
\end{example}

\noindent
{\sc Proof}:
Let $\theta_0$ be a support point of $G^*$.
There must be a positive constant $\delta$ such that
\[
f(x; G^*) 
\geq 
\delta \frac{(\theta_0)^x}{x!}  \exp( - \theta_0).
\]
Therefore, we have
\[
\EE^*\{  \log f(X; G^*) \}
\geq 
 \log (\theta_0)\, \EE^*(X) - \EE^* \{ \log (X!)\} + \log \delta- \theta_0.
\]
The condition $G^*(M) = 1$ easily leads to the finiteness of both $\EE^* (X)$
and $\EE^* \log (X!)$. Hence,
$\EE^*\{  \log f(X; G^*) \} > -\infty$.
Since the probability mass function is bounded from above,
we also have $\EE^*\{  \log f(X; G^*) \} < \infty$.
Therefore, $\EE^* | \log f(X; G^*) | < \infty$.

Since $f(x; B_\epsilon(G)) < 1$ for any $G$ and $\epsilon > 0$, 
\[
\EE^* [ \log \{ f(x; B_\epsilon(G)) / f(X; G^*) \}]
\leq
- \EE^*\{  \log f(X; G^*) \} < \infty.
\]
This is (KW3) for all $G \in \bar{\GG}$.
\qed \vs

Condition (KW3) is not always satisfied even when
the component distribution is Poisson.
Let $G(\{ A\})$ be the probability of $\theta \in A$ under
mixing distribution $G$.

\begin{example}
Let $G^*$ be a mixing distribution such that
\[
G^*( \{\log n\}) = c \{n (\log n) (\log\log n)^2\}^{-1}
\]
with some normalizing positive constant $c$, for $n=20, 21, \ldots$
Then under a Poisson mixture model, we have
\[
\EE^* \{\log f(X; G^*)\} = - \infty.
\]
Consequently, for any $\epsilon > 0$,
\[
\EE^*\{  \log f(X, B_\epsilon(1))/ f(X; G^*)\} =  \infty
\]
where $B_\epsilon(1)$ is the set of all mixing distributions
within $\epsilon$-distance of $f(x; 1)$.
%%% distance should be defined on mixing distribution space.
%%%   Need some correction here.

\end{example}

\noindent
{\bf Proof}:
The size of $c$ in this example does not affect the proof, 
so we take $c=1$. We choose this specific mixing distribution because
$\EE^* \{ \theta \} = \infty$.
In addition, it can be seen that
\ba
f(x; G^*)
&=&
\sum_{n=20}^\infty 
\big [
\{n (\log n) (\log\log n)^2\}^{-1} \frac{ (\log n)^x}{x!} \exp ( - \log n)
\big ]
\\
&=&
\frac{1}{x!} \sum_{n=20}^\infty \frac{(\log n)^{x-1}}{(n\log \log n)^2}.
\ea
It is seen that
\[
\sum_{n=20}^\infty \frac{(\log n)^{x-1}}{(n\log \log n)^2}
\approx
\int_{t=20}^\infty \frac{(\log t)^{x-1}}{(t \log \log t)^2}dt.
\]
The approximation is so precise that an error assessment
is possible but unnecessary. 
Changing the variable via $u = \log t$, we find
\[
\int_{t=20}^\infty \frac{(\log t)^{x-1}}{(t \log \log t)^2}dt
\leq
\int_{u=0}^\infty u^{x-1} \exp(-u) du
= (x-1)!
\]
The choice of $\log(20)$ as the lowest support point of
$G$ ensures that $\log(\log(20)) > 1$, which avoids some technicalities
in the above inequality.
Hence,
\[
\log f(x; G^*) \leq  \log \{ (x-1)!/x!\}  = -  \log x.
\]
Using $\EE^* = \EE_{G^*} \EE_\theta$, we find
\[
\EE^*\{  \log (X) \}
\geq
\EE_{G^*} \log \{ \EE_\theta X\}
=
\EE_{G^*} \{ \log \theta\}
= \infty.
\]
Therefore, 
\[
\EE^*\{  \log f(X; G^*) \} \leq - \EE^* \log (X) = - \infty.
\]

Let $\delta_x$ be a distribution with all the probability mass at $x$.
It can be seen that
\[
\dis(\delta_1, (1-\epsilon)\delta_1 + \epsilon \delta_x)
\leq
 \epsilon \int \exp( - \theta) d\theta = \epsilon.
 \]
 Therefore, $(1-\epsilon)\delta_1 + \epsilon \delta_x \in B_\epsilon(1)$
 for any $x$.
Hence, for any $x$ value,
\[
f(x; B_\epsilon(1)) \geq \epsilon f(x; x)
=
\frac{\epsilon x^x}{x!} \exp(- x) \approx \frac{\epsilon}{\sqrt{2\pi x}}
\]
by the Stirling formula. Because the Stirling formula is very accurate,
this implies that
\[
\log \{ f(x; B_\epsilon(1))/f(x; G^*)\} \geq \log (\epsilon/\sqrt{2\pi}) + (1/2) \log x.
\]
Hence, for any $\epsilon > 0$, we have
\[
\EE^* [ \log \{ f(X; B_\epsilon(1))/f(X; G^*)\}] 
\geq
\log (\epsilon/\sqrt{2\pi}) + (1/2)\EE^* \{\log X\}
= \infty.
\]
\qed

The point is that condition (KW3) on
a compact $\bar{\GG}$ places a severe
restriction on the KW proof.
Luckily, the Pfanzagl proof is free of this restriction.

\subsection{Consistency of nonparametric MLE under mixture model:
Pfanzagl approach}
\label{sub2.3}
The most demanding condition in the consistency proof of \cite{r1}
is (KW3), needed to validate the generalized Jensen's inequality \eqref{JensenG}.
Unfortunately, (KW3) under a mixture model is hard to verify,
as seen in the Poisson mixture example.
The Pfanzagl approach requires merely (W2) in comparison.
Here is the inequality in \cite{r2} that takes over the
role of Jensen's inequality.

\begin{lemma}
\label{lemmaF1}
Let $f(x)$ and $f^*(x)$ be density functions of any two distributions
with respect to some $\sigma$-finite measure $\mu$.
For any $u \in (0, 1)$, we have
\[
\EE^* \log \{ 1 + u [ f^*(X)/f(X) - 1]\} \geq 0,
\]
and equality holds if and only if $f^*(x) = f(x)$ almost surely with respect to the $f^*$ distribution.
\end{lemma}

\noindent
{\sc Proof}:
Let $Y =  \{ 1 + u [ f^*(X)/f(X) - 1]\} = (1-u) + u \{f^*(X)/f(X) \} $.
It can be seen that
\[
\log Y \geq (1-u) \log (1) +u \log \{f^*(X)/f(X)\} = u \log \{f^*(X)/f(X)\}.
\]
Hence,
\[
\EE^* \log (Y) \geq u \EE^*[ \log \{f^*(X)/f(X)\}] \geq 0,
\]
as required.
\qed
\vs

Now consider the mixture model $f(x; G)$ 
with its domain already extended to $\bar{\GG}$
as in the last subsection.

\vspace{1ex}
\begin{lemma}
\label{lemmaF2}
Assume that the mixture model is identifiable
and (KW2) is satisfied on $\bar{\GG}$.
Then, for any $G \in \bar{\GG}$ and $G \neq G^*$,  
there exists an $\epsilon>0$ such that
\be
\label{eqn20.42}
\EE^* \log \{ 1 + u [f(X; G^*)/f(X; B_\epsilon(G) )- 1]\}  > 0
\ee
where $G^*$ is the true mixing distribution.
\end{lemma}

\noindent
{\sc Proof}: 
Note that
\(
f(x; B_\epsilon(G))
\)
is a monotone increasing function of $\epsilon$, or it decreases to 
\(
f(x; G)
\)
for all $x$  as $\epsilon \to 0^+$.
In addition, given $u \in (0, 1)$,
\[
 \log \{ 1 + u [f(X; G^*)/f(X; B_\epsilon(G)) - 1]\}
 \geq \log (1 - u) > -\infty.
\]
Namely, this function has a finite lower bound,
which enables the use of Fatou's Lemma.
Hence,
\ba
& \lim_{\epsilon \to 0^+}
\EE^* \log \{ 1 + u [f(X; G^*)/f(X; B_\epsilon(G)) - 1]\}  \\
&
\vspace*{4ex}
\geq
\EE^* \log \{ 1 + u [f(X; G^*)/f(X; G) - 1]\}
> 0.
\ea
Consequently, there exists a positive $\epsilon$ value at which
\eqref{eqn20.42} holds.
\qed
\vs

This lemma is the counterpart of the generalized
Jensen's inequality in the KW approach. 
It is challenging to verify the validity of the
generalized Jensen's inequality under a general mixture model. 
In contrast, Lemma \ref{lemmaF2}
is valid provided the parameterization is continuous and identifiable.

\begin{thm}
\label{thm20.15}
Assume that the mixture model is identifiable
and (KW2) is satisfied on $\bar{\GG}$. We have
\[
\dis(\hat G_n, G^*) \to 0.
\]
%%% add word: almost surely?
\end{thm}

\noindent
{\sc Proof}: 
For any $\delta > 0$, let $B^c_\delta(G^*)$ be the distributions in ${\bar \GG}$
that are at least a distance $\delta$ from $G^*$.
Because ${\bar \GG}$ is compact, so is  $B^c_\delta(G^*)$.
Thus, the continuity condition implies that there exists a finite number
of $G_k, k=1, \ldots, J$, with corresponding $\epsilon_k$ such
that
\be
\label{eqn20.430}
B^c_\delta(G^*) \subset \cup_{k=1}^J B_k
\ee
where $B_k = \{G:  \dis(G, G_k) <  \epsilon_k\}$,
and
\be
\label{eqn20.43}
\EE^* \log \{ 1 + u [f(X; G^*)/f(X; B_k) - 1]\}  > 0.
\ee
By the strong law of large numbers, \eqref{eqn20.43} implies
\[
n^{-1} \sum_{i=1}^n \log \{ 1 + u [f(x_i; G^*)/f(x_i; B_k) - 1]\} > 0
\]
almost surely for $k=1, 2, \ldots, J$.
Consequently, we have
\ba
0 & < &
 \sum_{i=1}^n \log \{ 1 + u [f(x_i; G^*)/f(x_i; B_k) - 1]\}\\
& \leq&
\sup_{G \in B_k} \sum_{i=1}^n \log \{ 1 + u [f(x_i; G^*)/f(x_i; G) - 1]\}
\ea
almost surely for each $k=1, 2, \ldots, J$. Combining
this inequality with \eqref{eqn20.430}, we get
\[
0 < 
\sup_{G \not \in {B_\delta(G^*)} }
\sum_{i=1}^n \log \{ 1 + u [f(x_i; G^*)/f(x_i; G) - 1]\}
\]
almost surely.
By interpreting the summation in terms of the log-likelihood function, we find
\[
\ell_n( u G^* + (1-u) G) 
>
\ell_n(G)
\]
for all $G \in B^c_\delta(G^*)$ almost surely.
Since the likelihood function at any $G \not \in B_\delta(G^*)$
is smaller than the likelihood value at another mixing distribution 
$u G^* + (1-u) G$ that is a member of $\bar{\GG}$, 
the members of $B^c_\delta(G^*)$, all of which are at least a
$\delta$-distance away from $G^*$,
cannot possibly attain the supremum of $\ell_n(G)$.
Hence, the nonparametric MLE must reside in the
$\delta$-neighborhood of $G^*$ almost surely.
The arbitrarily small size of $\delta$ implies that $\dis(\hat G, G^*) \to 0$
almost surely as $n \to \infty$. This completes the proof.
\qed
\vs

In this proof, \cite{r2} took tactical advantage of the
linearity of the mixture model in mixing distributions:
\[
u f(x; G^*) + (1-u) f(x; G) = f(x; uG^*+(1-u) G),
\]
which is the density function of another mixture distribution. 

There is a limitation in the Pfanzagl result.
Consider the finite mixture model where $\GG$ is replaced by
$\GG_m$ for a given $m$.
The Pfanzagl result is no longer applicable because 
$uG^*+(1-u) G$ likely has more than $m$ support points even 
if both $G^*$ and $G$ have only $m$ support points.
In contrast, the KW proof leads to the
consistency of the MLE under finite mixture models
provided the corresponding conditions are satisfied.
In addition, the KW conditions under finite mixture models
are simple to verify, and they hold widely.

For finite mixture models, there
is another widely cited paper on the consistency of the MLE.

\subsection{Consistency of the MLE under finite mixture model:
Redner approach}
\label{sub2.4}

The generic results in \cite{r3} are not restricted to the
finite mixture model: the paper examines models lacking
full identifiability. Without identifiability, the parameter
estimator is inconsistent in general. However, in many situations,
the estimator may be regarded as consistent from a different angle.

Suppose a probability model
has its density function given by
\[
f(x; \theta_1, \theta_2) = \theta_1 \exp ( - \theta_1 x)
\]
for $x > 0$. Of course, this is simply an exponential distribution
with rate parameter $\theta_1$; parameter
$\theta_2$ is irrelevant. Given a set of \iid samples from a distribution
in this model, the MLE of $\theta_1$ is consistent, while there is no way
to have $\theta_2$ consistently estimated. At the same time,
the consistent estimation of
$\theta_2$ is unnecessary: it has no role in this population.
Let the distance between two vectors be
\[
\rho( (\theta_1, \theta_2), (\eta_1, \eta_2)) = | \theta_1 - \eta_1|.
\]
Then the MLE would satisfy, almost surely,
\[
\rho( (\hat \theta_1, \hat \theta_2), (\theta_1^*, \theta_2^*)) \to 0
\]
as the sample size $n \to \infty$.

Finite mixture models may not appear as trivial as in this example.
The density function of the two-component normal mixture in mean parameter
is given by
\[
f(x; \alpha, \theta_1, \theta_2) 
= (1-\alpha) \phi(x - \theta_1)+  \alpha \phi(x - \theta_2).
\]
The mixture with $\eta = (\alpha, \theta_1, \theta_2) = (0.3, 1, 2)$
is identical to the mixture with  $\eta = (0.7, 2, 1)$.
These two parameter vectors are apparently distinct in $\R^3$.
Hence, identifiability is lost when the model is parameterized
with this scheme. When $\eta$ is regarded as a
vector in $\R^3$ equipped with the Euclidean distance, the MLE
is not consistent.

Each $\eta$ in the above model has a corresponding
mixing distribution $G$. 
Let us define a distance for $\R^3$ (with the first component in [0, 1])
as
\[
\rho ( \eta_1, \eta_2) = \dis(G_1, G_2).
\]
Under this distance definition, we have
\[
\rho(\hat \eta, \eta^*) \to 0
\]
as $n \to \infty$ under appropriate conditions.

Since $\rho( \eta_1, \eta_2)=0$ does not lead to $\eta_1 = \eta_2$ 
in $\R^3$, $\rho(\cdot, \cdot)$ is not a mathematical distance. 
Regarding several distinct members of $\R^3$ as identical leads to
quotient topology. This turns out to be the core of \cite{r3}.

\subsection{Summary}
The consistency result of \cite{r2} seems perfect except for
the consistency of the MLE under finite mixture models.
The proof of \cite{r1} contains conditions that are not user-friendly.
However, these conditions become simple under finite mixture models.
Hence, the two papers perfectly complement each other.
\cite{r3} resorts to quotient topology to resolve nonidentifiability
and thereby provides another consistency proof for the finite mixture model.
\cite{r3}, however, assumes that $\Theta$ is a compact subset
of $\R^d$, making the result weaker. 

These papers do not consider only mixture models,
as this section may have suggested.
In this section we have substantially 
streamlined the conditions and conclusions and provided
additional insight in the context of mixture models.

\section{Consistency under finite normal mixture model}
One common requirement for the consistency of the MLE under a mixture model
is that $f(x; G)$ can be continuously extended to
$\bar{\GG}$. The first step of this extension is to have $f(x; \theta)$
continuously extended to include all $\theta$ on the boundary of $\Theta$.
This turns out to be impossible for the normal model with density
function
\[
\phi (x; \mu, \sigma) 
= \frac{1}{\sqrt{2\pi} \sigma}
\exp \left \{ - \frac{ (x - \mu)^2}{2 \sigma^2} \right \}.
\]
Because of this, none of the three approaches in the last section
is applicable to normal mixture models. 

The above issue is not the only obstacle.
The normal mixture model is not identifiable on $\GG$ unless
that space is reduced to $\GG_m$ for a prespecified $m$.
This section is devoted to the consistent estimation of $G$
under the finite normal mixture model.

\subsection{Finite normal mixture model with equal variances}
Consider the finite normal mixture model
where the component distributions share an equal but unknown variance:
\be
\label{EQnormal}
f(x; G; \sigma^2) =
\sum_{j=1}^m \alpha_j  \phi(x; \theta_j, \sigma^2).
\ee
Let  $G$ be the mixing distribution in component mean on $\Theta = \R$. 
The common variance $\sigma$ is structural with parameter space $\R^+$.

The log-likelihood function based on a set of \iid samples is given by
\[
\ell_n(G; \sigma^2) 
=
\sum_{i=1}^n \log f(x_i; G; \sigma^2).
\]
Here is a preliminary result similar to  but much strengthened
over that in \cite{r17}. We no longer confine the mean parameter
in a finite interval.

\begin{lemma}
\label{lemmaN1}
Let $(\hat G, \sigma^2)$ be a global maximum point of
the likelihood function $\ell_n(G; \sigma^2)$.
Then, there exist constants $0<\epsilon < \Delta < \infty$
such that as $n \to \infty$, the event sequence
$\{\epsilon \leq \hat{\sigma}^2 \leq \Delta\}$ occurs almost surely.
\end{lemma}

\noindent
{\sc Proof}:
It can be seen that $f(x; G; \sigma^2) \leq 1/\sigma$ for all $x$ and $G$.
When $\sigma^2 > \Delta$,
\[
\ell_n(G; \sigma^2) \leq - (1/2) n \log \Delta.
\]
Let $\bar x$ be the sample mean, and 
$s_n^2 = n^{-1} \sum_{i=1}^n (x_i - \bar x)^2$.
Then
\[
\ell_n(\bar x; s_n^2) \geq - n \log (s_n) -  (n/2).
\]
It can then be seen that
\ba
\ell_n(\bar x, s_n^2)
-
\ell_n(G; \sigma^2) 
&\geq&
\{ - n \log (s_n) -  (n/2) \} - \{ - (1/2) n \log \Delta\}\\
&=&
n \{ \log \Delta - \log (s_n) - (1/2) \}
\ea
uniformly for all $\sigma^2 > \Delta$.

Let $X$ be a random variable with
the true finite normal mixture distribution.
By the strong law of large numbers, 
$s_n^2$ almost surely converges to $\var(X)$
as $n \to \infty$.
Hence, when $\log \Delta > \log\{ \var(X)\} + (1/2)$, we have
\[
\ell_n(G; \sigma^2) < \ell_n(\bar x, s_n^2)
\]
almost surely for all $\sigma^2 > \Delta$.
This proves that the MLE for $\sigma^2$ is below this finite value $\Delta$
almost surely.

Next, because of the algebraic form of the normal density,
\be
\label{eqn45.11}
\log f(x; G, \sigma^2) \leq - \log (\sigma)
\ee
regardless of the actual value of $x$.
At the same time,
\[
f(x; G; \sigma^2) 
=
\sum_{j=1}^m \alpha_j  \phi(x; \theta_j, \sigma^2)
\leq
\max_j  \phi(x; \theta_j, \sigma^2).
\]
Hence, for any $G$ and $\sigma^2$, there is another upper
bound:
\[
\log f(x; G; \sigma^2) 
\leq  - \log \sigma -  (2 \sigma^2)^{-1} \min_{1 \leq j \leq m} (x - \theta_j)^2.
\]

Let $M$ be an arbitrary positive number and denote the
truncated $\theta$ value as
\[
\tilde \theta = 
\left \{
\begin{array}{cc}
-M & ~~~ \theta < -M;\\
\theta &~~~ |\theta| < M;\\
M &~~~ \theta > M.
\end{array}
\right .
\]
Let $\btheta= (\theta_1, \ldots, \theta_m)$ and
$\tilde \btheta= (\tilde \theta_1, \ldots, \tilde \theta_m)$.
The space of $\tilde \btheta$ is clearly compact given finite $M$.
For any $|x| \leq M$, we have
\be
\label{eqn45.12}
\log f(x; G, \sigma^2) 
\leq
- \log \sigma - (2 \sigma^2)^{-1} \min_{1 \leq j \leq m} (x - \tilde \theta_j)^2.
\ee

Applying \eqref{eqn45.11} for $|x_i| >M$ and \eqref{eqn45.12}
for $|x_i| \leq M$,
\be
\ell_n(G, \sigma^2)
\label{eqn45.13}
\leq
- n \log \sigma 
-
(2 \sigma^2)^{-1} \sum_{i=1}^n 
\big \{ \min_{1 \leq j \leq m} (x_i - \tilde \theta_j)^2 \big \} I( |x_i| \leq M).
\ee
Now we focus on the stochastic quantity
\[
h(\tilde \btheta, X)
=
\big \{ \min_{1 \leq j \leq m} (X - \tilde \theta_j)^2 \big \} I( |X| \leq M).
\]
Clearly, $h(\tilde \btheta, x) \leq M^2$ and it is equicontinuous in $\tilde \btheta$
for all $x$. 
Because the space of $\tilde \btheta$ is compact,
by the uniform strong law of large numbers of \cite{r16},
\[
n^{-1} \sum_{i=1}^n h(\tilde \btheta, x_i)
\to
\EE^* \{ h(\tilde \btheta, X)\}
\]
almost surely and uniformly in $\tilde \btheta$. 
Because $\EE^* \{ h(\tilde \btheta, X)\}$ is
smooth in $\tilde \btheta$ and it is clearly nonzero at each $\tilde \btheta$,
this implies
\[
\inf \EE^* \{ h(\tilde \btheta, X)\}  = \delta > 0
\]
where the infimum is over the compact space of $\tilde \btheta$.  
Applying this result to \eqref{eqn45.12}, we find
\[
\ell_n(G, \sigma) \leq - n \{ \log \sigma + \delta /\sigma^2 \}
\]
almost surely for all $\sigma$ and therefore
\be
\label{eqn45.14}
\ell_n(G, \sigma) -\ell_n(\bar x, s_n^2)
 \leq - n \{ \log \sigma  - \log (s_n) + \delta /\sigma^2 -1/2 \}
\ee
almost surely. 
When $\sigma^2$ is small enough, the upper bound
goes to negative infinity as $n \to \infty$. 
Hence, the maximum value of $\ell_n(G, \sigma)$ must be attained 
when $\sigma > \epsilon$ for some $\epsilon > 0$.
This completes the proof.
\qed
\vs

The key improvement of this lemma over that in \cite{r17} is
that here the parameter space of $\theta$ is the noncompact $\R$.

Based on this result, under
the finite normal mixture model with equal variance, the effective
component parameter space for $(\theta, \sigma^2)$
is $\R \times [\epsilon, \Delta]$
from the asymptotic point of view. Restricting the
space of $\sigma$ in this way leads to a compact parameter space.
On this component parameter space, we have
\[
\lim_{|\theta| \to \infty} f(x; \theta, \sigma^2) = 0.
\]
Hence, conditions (W2) and (W3) are satisfied after this restriction
and the KW approach is applicable.

\begin{thm}
\label{thmN1}
Under the finite normal mixture model \eqref{EQnormal}
with $m$ known, the MLE $(\hat G, \hat \sigma^2)$ is strongly consistent.
\end{thm}

The proof is simple.
Lemma \ref{lemmaN1} implies that the effective component parameter
space is $\R \times [\epsilon, \Delta]$.
The KW conditions on the reduced component parameter space
are satisfied. Hence, Theorem \ref{KWthm} can be applied to
give the consistency result.

\subsection{Finite normal mixture model with unequal variances}
The unequal-variance assumption does not exclude the possibility that
the true component variances are all equal. 
The density function is now given by
\be
\label{NEQnorm}
f(x; G) =
\sum_{j=1}^m \alpha_j  \phi(x; \theta_j, \sigma_j^2).
\ee
In this case, the mixing distribution $G$ is bivariate
on $\R \times \R^+$ and it mixes both mean and variance.
The log-likelihood has the same symbolic form:
\[
\ell_n(G) = \sum_{i=1}^n \log f(x_i; G).
\]

Consider the case where $m=2$. Let $(\theta_1, \sigma_1)=(0, 1)$,
$\alpha_1 = \alpha_2 = 0.5$. Let $\theta_2 = x_1$ and $\sigma_2 = 1/2k$
for $k=1, 2, \ldots$ 
Let $G_k$ be the corresponding mixing distribution. This
setup creates a sequence of mixing distributions $\{G_k\}_{k=1}^\infty$.
It can be seen that
\[
f(x_1; G_k) 
= 
\frac{0.5}{\sqrt{2\pi}}(2k) + \frac{0.5}{\sqrt{2\pi}} \exp( -\frac{ x_1^2}{2})
\geq \frac{k}{{2\pi}}
\]
and that for $i \geq 2$,
\ba
f(x_i; G_k) 
&=&
\frac{0.5}{\sqrt{2\pi}} (2k) \exp( - 2k^2(x_i-x_1)^2) 
+ \frac{0.5}{\sqrt{2\pi}} \exp( - \frac{x_i^2}{2})\\
&\geq&
\frac{1}{{2\pi}}\exp( - \frac{x_i^2}{2}).
\ea
Consequently, we have
\[
\ell_n(G_k) 
\geq
\log(k) - \frac{1}{2} \sum_{i=2}^n x_i^2 - n \log (2\pi).
\]
Clearly,  $\ell_n(G_k)  \to \infty$ as $k \to \infty$.
Hence, the limiting point of $G_k$ is one of the MLEs of $G$,
which is inconsistent.

There are some misconceptions in the literature. 
Since $\pr (X_1 = \theta) = 0$ for any given $\theta$ value,
one may suggest that the probability of having a degenerate MLE is zero.
This is false because $\theta = x_1$, in which $x_1$ is an
observed value, is no longer random after the observation.

In applications, EM-algorithm can be used to locate
many nondegenerate local maxima of $\ell_n(G)$.
The one with the largest likelihood value
is a locally consistent MLE of $G$ \citep{rr3,rr4}.
One may also use a consistent estimator, possibly
via method of moments, as an initial mixing distribution
for the EM-algorithm.
\cite{rr8} developed an approach to
test for global maximum, which can be useful
in the current context.
Experience shows that
such an estimator has good statistical properties, so
we should not write off this practice.

The inconsistency conclusion may not be an obstacle in many applications. 
Nonetheless, it is more satisfactory to have a foolproof method with
solid underlying theory that performs well in applications.
In the literature there are two approaches to consistent estimation
based on likelihood.
One is the constraint MLE proposed by \cite{r9}. Simply put,
it reduces the component parameter space of $\sigma$.
The result of the last subsection may be regarded as its
simplest case.

This subsection focuses on the penalty method applied to $\sigma$. 
The penalty is also a prior on $\sigma$ or a regularization measure.
The first largely successful proof of consistency
for the penalized MLE is in \cite{r10}, following its
proposal in \cite{r12}.
\cite{r6} provided a successful complete proof that
is simplified and improved here.

\subsubsection{Penalized likelihood.}

The inconsistency of the MLE under the finite normal mixture
model is largely due to nonregularity. 
Hence, regularizing the likelihood is a natural way
to gain consistency of the altered MLE.
The regularization is itself in the form of a penalized
likelihood as follows:
\[
\tilde \ell_n(G) = \ell_n(G) + p_n(G).
\]
The mixing distribution is then estimated
 by one of the global maxima of  $\tilde \ell_n(G)$ over $\GG_m$:
\[
\tilde G = \arg \max \tilde \ell_n(G).
\]
The uniqueness is a natural consequence of the
subsequent discussion.
We denote the component means in $\tilde G$
as $\tilde \theta_j$ and so on. 

\subsubsection{Technical lemmas.}
The following lemma
provides a technical basis for the size of the penalty.

\begin{lemma}
\label{lemmaN2}
Let $x_1, \ldots, x_n$ be a set of $n$ \iid observations from an
absolute continuous distribution $F$
with density function $f(x)$. Assume that $f(x)$  is continuous
and $M = \sup_x f(x) < \infty$. 
Let $F_n(x) = n^{-1} \sum_{i=1}^n I(x_i \leq x)$
be the empirical distribution function.

Then, as $n\to \infty$ and almost surely, for any given 
$\epsilon > 0$, 
\[
\sup_{\theta \in \R} \{ F_n(\theta + \epsilon) - F_n(\theta)\}
\leq 2M \epsilon + 8 n^{-1}\log n.
\]
\end{lemma}

\noindent
{\sc Proof}:
Since $F(x)$ is continuous, there exist
$\eta_0, \eta_1, \ldots, \eta_n$ such that $F(\eta_i) = j/n$ for $0 < j < n$
with $\eta_0 = - \infty$ and $\eta_n = \infty$. This ensures that
for any $\theta$ value, there exists a $j$ such that
$\eta_{j-1} < \theta \leq \eta_j$.
Therefore,
\begin{eqnarray}
\sup_\theta \{ F_n(\theta + \epsilon) - F_n(\theta)\}
&\leq&
\max_j \{ F_n(\eta_j + \epsilon) - F_n(\eta_{j-1})\}
\nonumber \\
&\leq&
\max_j |\{ F_n(\eta_j + \epsilon) - F_n(\eta_{j-1})\} -  \{ F(\eta_j + \epsilon) - F(\eta_{j-1})\}|
\nonumber \\
&& 
+ \max_j  \{ F(\eta_j + \epsilon) - F(\eta_{j-1})\}.
\label{eqn45.10}
\end{eqnarray}
The task of the proof is to find appropriate bounds for these two terms.
First,
\[
 F(\eta_j + \epsilon) - F(\eta_{j-1})
 \leq
 \{ F(\eta_j + \epsilon) - F(\eta_{j})\} +  \{ F(\eta_j) - F(\eta_{j-1})\}.
\]
Since $F(\eta_j) - F(\eta_{j-1}) = n^{-1}$
and $F(\eta_j + \epsilon) - F(\eta_{j}) \leq M\epsilon$ by the mean value
theorem, we have
\[
\max_j  \{ F(\eta_j + \epsilon) - F(\eta_{j-1})\} \leq M\epsilon + n^{-1}.
\]

Let
\(
Y_i = I( \eta_{j-1} < x_i \leq \eta_j + \epsilon)
\)
and write $\Delta_j = n^{-1} \sum_i \{Y_i - \EE Y_i\}$.
The first term in (\ref{eqn45.10}) equals $\max_j \Delta_j$.
By applying Bernstein's inequality to $\Delta_j$ followed
by the Borel--Cantelli Lemma as in \cite{r13}, we
get
\[
\max_j \Delta_j < M\epsilon + 8 n^{-1} \log n.
\] 
Combining the two bounds leads to
\be
\label{eqn45.101}
\sup_\theta \{ F_n(\theta + \epsilon) - F_n(\theta)\}
\leq 2M\epsilon + n^{-1} +  8 n^{-1} \log n
\ee
almost surely. Because $n^{-1}$ is a high-order term
compared to $n^{-1} \log n$,
it is absorbed into the latter.
\qed

\vs
The proof remains solid when $\epsilon$ depends on $n$. 
Because of this, $\epsilon$ is allowed to take an arbitrarily small value
without invalidating the inequality.
Technically, Lemma \ref{lemmaN2}
leaves a zero-probability event for each value of $\epsilon$
on which the upper bound is violated.
The union of these zero-probability events over $\epsilon$
does not have to be a zero-probability event. 
However, since $\sup_\theta \{ F_n(\theta + \epsilon) - F_n(\theta)\}$
is monotone in $\epsilon$, this technicality is easily resolved.

\begin{lemma}
\label{lemma3.3}
The upper bound in Lemma \ref{lemmaN2} 
after a minor alteration,
\[
\sup_{\theta \in \R} \{ F_n(\theta + \epsilon) - F_n(\theta)\}
\leq 2M \epsilon + 10 n^{-1}\log n,
\]
holds uniformly for all $\epsilon > 0$ almost surely.
\end{lemma}

This lemma shows that \iid observations from a population
with a bounded density function spread out evenly almost surely.
This result extends a corresponding result in
\citep{r6} to cover generic distribution $F$.

\subsubsection{Choice of penalty.}
Imposing the following three properties on the penalty function makes
the penalized MLE consistent:

\vs
\noindent
{\bf P1}. Additivity: $p_n(G) = \sum_{j=1}^m \tilde p_n(\sigma_j)$.

\vs
\noindent
{\bf P2.} Uniform upper bound: $\sup_{\sigma > 0} [ \tilde p_n(\sigma)]^+ = o(n)$; 
individual lower bound: $\tilde p_n(\sigma) = o(n)$ for each $\sigma > 0$.

\vs
\noindent
{\bf P3.} Sufficiently severe: $\tilde p_n(\sigma) < (\log n)^2 \log (\sigma)$
for $\sigma < n^{-1}\log n$ when $n$ is large enough.

\vspace{1ex}
The first property allows for a simple discussion and
straightforward numerical solution to the penalized MLE.
%%%
The upper and lower bounds in {\bf P2} 
prevent the likelihood from being seriously inflated or deflated
at any $G$.
Property {\bf P3} requires the size of the penalty to be
large enough to prevent $\sigma_j \approx 0$ in the penalized
MLE of $G$.
One possible $\tilde p_n$ is
\[
\tilde p_n(\sigma) = - n^{-1} \{ \sigma^{-2} + \log \sigma^2\}.
\]
This penalty function goes to negative infinity when 
$\sigma \to 0$ or $\sigma \to \infty$.
It is minimized when $\sigma = 1$. 
The upper-bound condition in ${\bf P2}$
on $[\tilde p_n(\sigma)]^+$ is certainly satisfied. 
The lower-bound condition in ${\bf P2}$
is clearly satisfied. 
As for {\bf P3},  when $\sigma \to 0^+$, $\sigma^{-2} \to \infty$ much
faster than $|\log \sigma|$.  Hence, the penalty is much more severe
than the required $4 (\log^2 n) \log \sigma$.
In applications, $\sigma$ should be replaced by $\sigma/s_n$, where
$s_n$ is the sample variance, to retain scale invariance.
This penalty function also represents a prior Gamma distribution placed on $\sigma^{-2}$.
This form of penalty is very convenient for the EM-algorithm
popularly used for numerical computation and on Bayes analysis
as in \cite{rr5,r12}.

\subsubsection{Consistency of the penalized MLE.}
We now outline the proof for $m=2$; the
general case is similar and omitted.
Let $K^* = \EE^* \{ \log f(X; G^*)\}$ and $M = \sup_x f(x; G^*)$. 
Select a sufficiently small $\epsilon_0$ such that

\vs\vs
\noindent
(1) $4M \epsilon_0 (\log \epsilon_0)^2 \leq 1$;
(2) $ (1/2)(\log \epsilon_0)^2  + \log (\epsilon_0) \geq 4 - 2 K^*$.
\vs\vs

Without loss of generality, we assume $\sigma_1 \leq \sigma_2$.
Partition the mixing distribution space $\GG_2$ into
\bi
\item[]
$
\Gamma_1 = \{G: \sigma_1 \leq \sigma_2 \leq \epsilon_0\}$;
\item[]
$
\Gamma_2 = \{ G: \sigma_1 < \tau_0, \sigma_2 > \epsilon_0\}
$
for some constant $\tau_0 < \epsilon_0$ to be specified;
\item[]
$\Gamma_3= \{ \Gamma_1 \cup \Gamma_2\}^c$.
\ei

The overall strategy is to show that the penalized MLE
is not in $\Gamma_1$ or $\Gamma_2$ but in $\Gamma_3$.
The final conclusion follows from the
fact that the finite normal mixture model on $\Gamma_3$ fits into
the KW proof.

\vs\vs\noindent
{\bf Step I}. The following lemma says it all. 

\begin{lemma}
\label{lemmaN4}
In the current setting, 
\(
\sup_{G \in \Gamma_1}  \tilde \ell_n(G) - \ell_n(G^*)
\to - \infty
\)
almost surely as $n \to \infty$.
\end{lemma}

\noindent
{\sc Proof}: Define
\(
A_j  =  \{i: |x_i-\theta_j| < |\sigma_j \log \sigma_j | \}
\) for $j=1, 2$.
Partition the entries in $\ell_n(G)$ to get
\[
\ell_n(G) = \ell_n(G; A_1) + \ell_n(G; A_1^cA_2) +\ell_n(G; A_1^cA_2^c)
\]
where $\ell_n(G; A) = \sum_{i \in A} \log f(x_i; G).$

Denote the number of observations in set $A$ as $n(A)$. 
Since the mixture density function is bounded by $\sigma_1^{-1}$
for $G \in \Gamma_1$,
\[
\ell_n( G ; A_1) \leq  - n(A_1) \log(\sigma_1).
\]
Applying Lemma \ref{lemma3.3} with 
$\epsilon = \sigma_1 \log (1/\sigma_1)$,
we get 
\be
\label{eq99}
n(A_1) \leq - 2nM \sigma_1 \log (\sigma_1) + 10 \log n.
\ee
Hence,
\begin{equation}
\label{eqn45.36}
\ell_n(G; A_1) 
\leq
2 nM \sigma_1 (\log \sigma_1)^2 - 10 (\log n)  \log (\sigma_1).
\end{equation}
By {\bf P3},  $\tilde p_n(\sigma_1) < (\log n)^2 \log (\sigma_1)$.
Hence, (\ref{eqn45.36}) leads to
\ba
\ell_n(G; A_1) + \tilde{p}_n(\sigma_1)
&\leq &
2nM \sigma_1 (\log \sigma_1)^2 -  \{10(\log n)-  (\log n)^2\} \log ( \sigma_1)
\\
&\leq &
2nM \sigma_1 (\log \sigma_1)^2.
\ea
Similarly,
\[
\ell_n(G; A_1^cA_2) + \tilde{p}_n(\sigma_2)
\leq
2nM \sigma_2 (\log \sigma_2)^2.
\]

Let $\phi(\cdot)$ be the density function of the standard normal.
The  observations falling outside both $A_1$ and $A_2$ have
log-likelihood contributions that are bounded by
\[
\log \{ \frac{\alpha_1}{ \sigma_1} \phi( - \log \sigma_1) + 
\frac{\alpha_2}{\sigma_2} \phi( - \log \sigma_2)\} 
\leq
- \log (\epsilon_0) -  \frac{1}{2}(\log \epsilon_0 )^2.
\]
At the same time, by (3.9) and for small enough $\epsilon_0$ and sufficiently
large $n$,
\[
n( A_1^cA_2^c) \geq n - \{ n( A_1) + n( A_2)\}  \geq  \frac{n}{2}.
\]
Hence, we get the third bound:
\[
\ell_n(G; A_1^cA_2^c)
 \leq - \frac{n}{2} \{ \log (\epsilon_0) 
+ \frac{1}{2} (\log \epsilon_0)^2\}.
\]

Combining the three bounds and remembering how $\epsilon_0$
was selected, we conclude that when $G \in \Gamma_1$,
\begin{eqnarray*}
\tilde \ell_n(G) 
&=& 
\{ \ell_n(G; A_1) + \tilde{p}_n(\sigma_1) \}
	+ \{ \ell_n(G; A_1^cA_2) - \tilde{p}_n(\sigma_2)\}
	+ \ell_n(G;A_1^cA_2^c) \\
& \leq & 
4 M n \epsilon_0 (\log \epsilon_0)^2
  -  \frac{n}{2} \{ \frac{1}{2}(\log \epsilon_0)^2 + \log (\epsilon_0)\} \\
& \leq &
n -  \frac{n}{2} (4 - 2K^*) \\
& = & 
 n (K^*-1).
\end{eqnarray*}
The last few inequalities hold by the tactical choice of $\epsilon_0$.

By the strong law of large numbers, 
$n^{-1} \tilde \ell_n (G^*) \to K^*$ almost surely. 
The last inequality is then simplified to 
\[
\sup_{G \in \Gamma_1} \tilde \ell_n(G) - \tilde \ell_n( G^*) 
\leq  -n \to -\infty.
\]
This completes the proof.
\qed\vs

Remark: I have omitted ``almost surely'' in the proof for ease
of presentation.
\vs

\noindent
{\bf Step II:}
The penalized MLE 
of $G$ is almost surely not inside $\Gamma_2$,
for an appropriately chosen $\tau_0$.
The choice may depend on $G^*$ but not on the sample
size $n$.
Let $\bar{\Gamma}_2$ be a compactified $\Gamma_2$ allowing
$\sigma_1 = 0$ and $\alpha_1+\alpha_2 < 1$.
Define, for any $G \in \bar{\Gamma}_2$,
\be
\label{eqn45.gg}
g(x; G)
= 
\alpha_1 \phi(x; \theta_1, 2\epsilon_0^2)
+ 
\alpha_2 \phi(x; \theta_2, \sigma_2^2).
\ee
On $\bar{\Gamma}_2$, $\sigma_2$ has a nonzero lower bound. 
Thus, $g(x; G)$ is bounded although $\sigma_1 =0$ is allowed.

Without loss of generality, $\tau_0$ is
small enough such that the true mixing
distribution $G^* \not \in \Gamma_2$. 
%%% some thinking needed here.
Hence, applying Jensen's inequality, we also have
\[
\EE^* \log \{g(X; G)/f(X; G^*)\} < 0
\]
for all $G \in \bar{\Gamma}_2$.
Using a slightly different symbol from $\ell$, we define
\[
l_n(G) = \sum_{i=1}^n \log \{ g(x_i; G)\}
\]
on $\bar{\Gamma}_2$.
By the strong law of large numbers
and the newly established Jensen's inequality,
\[
n^{-1} \{ l_n(G) - \ell_n(G^*) \} \to \EE^* \log \{g(X; G)/f(X; G^*)\} < 0.
\]
Further exploring this key conclusion leads to the following lemma.

\begin{lemma}
\label{lemma3.5}
Consider a set of $n$ \iid observations from $f(x; G^*)$ and
the function $g(x; G)$ defined by \eqref{eqn45.gg}.
Let $l_n(G) = \sum_{i=1}^n \log \{ g(x_i; G)\}$.
We then have
\begin{equation}
\label{eqn45.31}
l_n(G) - \ell_n(G^*) \leq - n \delta(\epsilon_0)
\end{equation}
for some $\delta(\epsilon_0) > 0$ almost surely.
\end{lemma}

\noindent
{\bf Proof}: 
Let
\[
g(x; G, \epsilon) 
= 
\sup \{g(x; G'):  G' \in \bar{\Gamma}_2,~ \dis(G', G) < \epsilon\}.
\]
Because $\sigma_2 > \tau_0 > 0$, we still have
\(
g(x; G, \epsilon) \leq  1 + \tau_0^{-1}.
\)
Hence, $\EE^* \log g(X; G, \epsilon)< \infty$.

Clearly, $\EE^* \log g(X; G, \epsilon) > -\infty$.
Therefore,
\(
\EE^*\log  \{ g(X; G, \epsilon)/f(X; G^*)\}
\)
is well defined.
Let $\epsilon \to 0^+$; by the monotone
convergence theorem we find
\[
\lim_{\epsilon \downarrow 0} 
\EE^*\log \{ g(X; G, \epsilon)/f(X; G^*)\}
\leq 
\EE^* \log  \{ g(X; G)/f(X; G^*)\}
< 0.
\]
Next, note that $\bar{\Gamma}_2$ is
compact based on the distance $\dis(\cdot, \cdot)$.
There is a finite number of $G$ and $\epsilon$ such
that
\[
\bar{\Gamma}_2 \subset \cup_{j=1}^J \{ G: \dis(G, G_j) \leq \epsilon_j\}
\]
and for each $j=1, \ldots, J$,
\[
\EE^* \{ g(X; G_j, \epsilon_j)/f(X; G^*)\} < 0.
\]
This leads to the claim of this lemma:
\[
l_n(G) - \ell_n(G^*) \leq - n \delta(\epsilon_0)
\]
for some $\delta(\epsilon_0)>0$ whose size depends on 
the size of $\epsilon_0$.
\qed \vs

Let us connect $l_n(G)$ to $\ell_n(G)$ on the space 
$\bar {\Gamma}_2$
and refine this result to obtain the major result of
this step.

\begin{lemma} 
\label{lemma}
As $n \to \infty$,
\(
\sup_{G \in \Gamma_2} \tilde{\ell}_n(G) - \tilde{\ell}_n(G^*) \to -\infty.
\)
\end{lemma}

\noindent
{\bf Proof}: 
Retain the definition
$A_1 = \{ i: |x_i - \theta_1| \leq \sigma_1 \log (1/\sigma_1) \}$.
For each $i \in A_1$, we have
\[
f(x_i; G) \leq  (1/\sigma_1) g(x_i; G).
\]
Therefore, the log-likelihood contribution of each observation in $A_1$ 
is
\[
\log \{ f(x_i; G) \} \leq \log (1/\sigma_1) + \log\{ g(x_i; G)\}.
\]
The observed values of the observations not in $A_1$
satisfy $|x - \theta_1| \geq |\sigma_1 \log \sigma_1|$.
For these $x$ values, we have
\[
\frac{(x- \theta_1)^2}{2 \sigma_1^2}
\geq
\frac{(x- \theta_1)^2}{4 \sigma_1^2}
+ \frac{1}{4} (\log \sigma_1)^2
\geq
\frac{(x- \theta_1)^2}{4 \epsilon_0^2}
+ \frac{1}{4} (\log \sigma_1)^2.
\]
Consequently,
\bea
\frac{1}{\sigma_1} 
\exp \big \{ - \frac{(x- \theta_1)^2}{2 \sigma_1^2} \big \} 
&\leq &
\exp \big \{ - \frac{(x- \theta_1)^2}{4 \epsilon_0^2} \big \} 
\times 
\exp \big \{ - \frac{1}{4} (\log \sigma_1)^2 - \log \sigma_1 \big \}
\nonumber \\
&= &
\exp \big \{ - \frac{(x- \theta_1)^2}{4 \epsilon_0^2} \big \} 
\times 
\exp \big \{ - \frac{1}{4} (\log \sigma_1+2)^2 +1 \big \}.
\label{eq100}
\eea
The factor $1/\sigma_1$ has been turned
into $\exp ( - \log \sigma_1)$ in the first line of the
above derivation.
For a small enough $\epsilon_0$,
\[
\exp \big \{ - \frac{1}{4} (\log \sigma_1+2)^2 +1 \big \} \leq \frac{1}{2 \epsilon_0}
\]
when $\sigma_1 \leq \epsilon_0$.
Hence, \eqref{eq100} leads to, for those $x$ values not in the set $A_1$,
\[
\phi(x; \theta_1, \sigma_1^2)
\leq
\phi(x; \theta_1, 2 \epsilon_0^2)
\]
and therefore,
\ba
f(x; G) 
&=&
\alpha_1 \phi(x; \theta_1, \sigma_1^2) + \alpha_2 \phi(x; \theta_2, \sigma_2^2)\\
&\leq&
\alpha_1 \phi(x; \theta_1, 2 \epsilon_0^2) + \alpha_2 \phi(x; \theta_2, \sigma_2^2)\\
&=&
g(x; G).
\ea
In summary, when $i \not \in A_1$, its log-likelihood contributions
\[
\log f(x_i; G) \leq \log \{g(x_i; G)\}. 
\]

Combining the cases for the observations in and not in $A_1$, we find
\[
\ell_n(G) \leq n(A_1) \log (1/\sigma_1) + \sum \log \{g(x_i; G)\}. 
\]
This leads to
\[
\sup_{G \in \Gamma_2} \tilde \ell_n(G) 
\leq
\sup_{G \in \Gamma_2} \{l_n(G) + \tilde p_n(\sigma_2)\}
+
\sup_{G \in \Gamma_2} \{ n(A_1) \log (1/\sigma_1) + \tilde p_n(\sigma_1)\}.
\]
Reusing the bound \eqref{eq99} on $n(A_1)$ together with {\bf P3},
we get
\[
\sup_{G \in \Gamma_2} \{ n(A_1) \log (1/\sigma_1) + \tilde p_n(\sigma_1)\}
<  2 M n \tau_0 (\log \tau_0)^2.
\]
Hence, the proof of the lemma is reduced to showing that
\be
\label{eqn45.32}
\big [ \sup_{\Gamma_2} \{l_n(G) + \tilde p_n(\sigma_2)\} 
+ 
2 M n \tau_0 (\log \tau_0)^2 \big ]
- \tilde \ell_n(G^*) 
< 0.
\ee
Because $[p_n(\sigma_2)]^+ = o_p(n)$ by choice,
it suffices to show that
\begin{equation}
\label{eqn45.33}
\sup_{\Gamma_2} l_n(G) - \ell_n(G^*) 
\leq - \delta n
\end{equation}
for some $\delta > 2 M \tau_0 (\log \tau_0)^2$.
This is implied by Lemma \ref{lemma3.5} when a sufficiently small $\tau_0$
is chosen, after the choice of $\epsilon_0$.
\qed\vs

The $g(x; G)$ used here is more convenient than that in \cite{r6}.
The proofs so far have successfully
excluded the possibility that the penalized MLE of $G$
falls in $\Gamma_1\cup \Gamma_2$.
Finishing the consistency proof is a simple task.

\begin{thm}
\label{thm45.33}
The penalized MLE of $G$ is consistent:
$\tilde G \to G^*$ almost surely as $n \to \infty$.
\end{thm}

Once the mixing distribution is restricted to $\Gamma_3$,
the KW conditions are satisfied. Hence, the restricted MLE is consistent.
On this space, the penalty is of size $o_p(n)$. Hence, the
penalized MLE remains consistent.
Note that the KW conditions are easy
to verify for finite normal mixture models on $\Gamma_3$.

\subsection{Summary}
The finite normal mixture model does not satisfy the KW
conditions, and the MLE as defined by \eqref{eqn20.mle} together
with subsequent remarks is inconsistent.
The consistency of the penalized MLE has only recently been
solidly proved.  Yet these facts are often overlooked;
we have made these results more accessible to researchers in various disciplines.
The penalized MLE under a multivariate
normal mixture has also been shown to be consistent by \cite{r18}
with a minor correction by \cite{r19}. Note that
the proof of the multivariate case can be substantially simplified
using the new techniques in this paper. 
Finally, being consistent is a minimum requirement
in statistical data analysis. 
The proper estimation of the mixing distribution under a finite mixture model
requires a very large sample size when the subpopulations
are not well separated, as indicated by the simulation study of \cite{r18}.

\section{Concluding remarks}
We have discussed the consistency of the MLE under mixture models given \iid
observations. 

When there are no
restrictions on the mixing distribution, the nonparametric
MLE is consistent under the minimum identifiability condition
(KW1) and the continuity condition (KW2) based on 
Pfanzagl's proof. In addition, (KW2) is implied by (W2) and (W4).

The conclusion via the KW proof is most
useful when applied to finite mixture models.
In this case, the KW conditions are implied by Wald conditions
(W2)--(W4) and (KW1). The MLE of the mixing distribution
under a finite mixture is consistent.

Under the finite normal mixture model with equal variances,
the MLE is consistent. Under the finite normal mixture model in both mean and variance,
the MLE is not consistent. When a penalty satisfying {\bf P1}--{\bf P3}
is applied to the log-likelihood function, the penalized MLE
is consistent.

It is curious that the consistency proof given by Pfanzagl
for the nonparametric MLE under a mixture model requires
a somewhat weaker set of conditions than its parametric
counterpart given by Wald: (KW1), (W2), and (W4), but not (W3). 
Because Pfanzagl's result is not applicable to
finite mixture models, it does not cover the
result of Wald or that of KW.
The KW proof makes that of Wald a special case.

This paper is novel at streamlining the conditions,
conclusions and the proofs related to the consistency of 
the MLE under mixture models.
The proofs of \cite{r1,r2,r3,r4} are substantially simplified
by connecting all of them with Theorem \ref{Wald}.
The proofs focus on essential ideas and leaves complex conditions
out for separate discussions.
The Wald consistency result is strengthened by
de-requiring $E^*|\log f(X; \theta^*)| < \infty$.
The KW consistency conclusion is found most useful to
finite mixture models. Its generic conclusion is less useful
because its (KW3) condition is difficult to verify
or not satisfied as testified by a Poisson mixture example.
The paper follows the existing line of proofs for
the consistency of the MLE under finite normal mixture models.
When the variance is a structural parameter, 
the conclusion is strengthened by de-require the space of the 
mean parameter being compact. 
The consistency proof of the penalized MLE under
finite normal mixture models is simplified though
developing a more generic concentration inequality
in Lemma \ref{lemma3.3} and introducing a
more convenient function $g(x; G, \epsilon)$ in Lemma \ref{lemma3.5}.

%\section*{Acknowledgements}

\end{document}